\documentclass[11pt]{article}
\usepackage{latexsym,amssymb}
\usepackage{amsmath}
\usepackage[mathscr]{eucal}

\newtheorem{Theorem}{\bf Theorem}[section]
\newtheorem{Lemma}{\bf Lemma}[section]
\newtheorem{Proposition}{\bf Proposition}[section]  
\newtheorem{Corollary}{\bf Corollary}[section]
\newtheorem{Remark}{\bf Remark}[section]

\newtheorem{Definition}{\bf Definition}[section]
\newtheorem{Condition}{\bf Condition}[section]

\newenvironment{theorem}{\begin{Theorem}$\!\!\!$}{\end{Theorem}}
\newenvironment{lemma}{\begin{Lemma}$\!\!\!$}{\end{Lemma}}
\newenvironment{proposition}{\begin{Proposition}$\!\!\!$}{\end{Proposition}}
\newenvironment{corollary}{\begin{Corollary}$\!\!\!$}{\end{Corollary}}
\newenvironment{remark}{\begin{Remark}$\!\!\!$}{\end{Remark}}
\newenvironment{definition}{\begin{Definition}$\!\!\!$}{\end{Definition}}

\numberwithin{equation}{section}
\begin{document}
\title{Power concavity for elliptic and parabolic
boundary value problems on rotationally symmetric domains}
\author{Kazuhiro Ishige, Paolo Salani and Asuka Takatsu}

\date{}
\maketitle
\begin{abstract}
We study power concavity of rotationally symmetric solutions to 
elliptic and parabolic boundary value problems on rotationally symmetric domains 
in Riemannian manifolds. 
As applications of our results to the hyperbolic space ${\bf H}^N$ 
we have: 
\begin{itemize}
  \item
  The first Dirichlet eigenfunction on a ball in ${\bf H}^N$ is strictly positive power concave;
  \item
  Let $\Gamma$ be the heat kernel on ${\bf H}^N$. Then 
  $\Gamma(\cdot,y,t)$ is strictly log-concave on ${\bf H}^N$ for $y\in {\bf H}^N$ and $t>0$. 
\end{itemize}
\end{abstract}
\vspace{40pt}
\noindent Addresses:

\smallskip
\noindent 
K. I.: Graduate School of Mathematical Sciences, The University of Tokyo, 3-8-1 
Komaba, Meguro-ku, Tokyo 153-8914, Japan\\
\noindent 
E-mail: {\tt ishige@ms.u-tokyo.ac.jp}\\

\smallskip
\noindent 
P. S.: Dipartimento di Matematica ``U. Dini'', 
Universit\`a di Firenze, viale Morgagni 67/A, 50134 Firenze\\
\noindent 
E-mail: {\tt paolo.salani@unifi.it}\\

\smallskip
\noindent 
A. T.: Department of Mathematical Sciences,  Tokyo Metropolitan University, 
Minami-osawa, Hachioji-shi, Tokyo 192-0397, Japan\\
\noindent 
E-mail: {\tt asuka@tmu.ac.jp}\\
\newpage
%%%%%%%%%%%%%%%%%%%%%%%%%%%%%%%%%%%%
%%%%%%%%%%%%%%%%%%%%%%%%%%%%%%%%%%%%
\section{Introduction}
%%%%%%%%%%%%%%%%%%%%%%%%%%%%%%%%%%%%
%%%%%%%%%%%%%%%%%%%%%%%%%%%%%%%%%%%%

Concavity of solutions to elliptic and parabolic boundary value problems on convex domains in Euclidean space 
is a classical subject and has fascinated many mathematicians. 
%%%
The literature is large
and we just refer to the classical monograph by Kawohl~\cite{KawohlBook} 
and the papers 
\cite{BL}, \cite{DD}, \cite{GK}, 
\cite{ILS}--\cite{Kawohl}, \cite{Kennington}, \cite{Korevaar},  \cite{LV} and
%, \cite{SWW} and \cite{Shih}, 
some of which are closely related to this paper and the others include recent developments in this subject. 
In this regard, let us recall some results on power concavity properties of solutions 
to elliptic and parabolic boundary value problems.
\begin{itemize}
  \item[(a)] 
  Let $\Omega$ be a bounded convex domain in ${\bf R}^N$. 
  Let $u\in C^2(\Omega)\cap C(\overline{\Omega})$ be a solution to 
  \[
  \left\{\begin{array}{ll}
  -\Delta  u=\lambda u^\gamma\quad&\mbox{in}\quad \Omega,\vspace{5pt}\\
  u>0\quad&\mbox{in}\quad \Omega,\vspace{5pt}\\
  u=0\quad&\mbox{on}\quad\partial \Omega,
  \end{array}\right.
  \]
  where $\lambda>0$ and $0\le\gamma\leq 1$.
  Then $u$ is $\alpha$-concave in $\Omega$ with $\alpha=(1-\gamma)/2$.  
  See e.g. \cite{Kawohl0, Kawohl, Kennington}.
  In particular, 
  if $\phi$ is \emph{the first Dirichlet eigenfunction} for $-\Delta$ on $\Omega$,
  namely $\phi$ satisfies 
  \[
  \left\{\begin{array}{ll}
  -\Delta  \phi=\lambda_1(\Omega)\phi\quad&\mbox{in}\quad \Omega,\vspace{5pt}\\
  \phi>0\quad&\mbox{in}\quad \Omega,\vspace{5pt}\\
  \phi=0\quad&\mbox{on}\quad\partial \Omega,
  \end{array}\right.
  \]
  then $\phi$ is $0$-concave (i.e. log-concave) in $\Omega$. 
  Here $\lambda_1(\Omega)$ is  the \emph{first Dirichlet eigenvalue} for $-\Delta$ on $\Omega$. 
  See e.g. \cite{BL, Korevaar}.
  \item[(b)] 
  Let $\Omega$ be a convex domain in ${\bf R}^N$. 
  Let 
  \[
  u\in C^{2,1}(\Omega\times(0,\infty))\cap C(\overline{\Omega}\times(0,\infty))
  \cap BC(\Omega\times[0,\infty))
  \]
  be a nonnegative  solution to
  \[
  \left\{
  \begin{array}{ll}
  \partial_t u=\Delta u-\lambda u^\nu & \quad\mbox{in}\quad\Omega\times(0,\infty),\vspace{3pt}\\
  u=0 & \quad\mbox{on}\quad\partial\Omega\times(0,\infty)\quad\mbox{if}\quad\partial\Omega\not=\emptyset,\vspace{3pt}\\
  u(x,0)=\varphi(x) & \quad\mbox{in}\quad\Omega,
  \end{array}
  \right.
  \]
  where $\lambda\ge 0$, $\nu\ge 1$ and $\varphi$ is a bounded, continuous and nonnegative function in $\Omega$.
  Then $u(\cdot,t)$ is log-concave in $\Omega$ for $t>0$ provided that $\varphi$ is log-concave in $\Omega$. 
  See e.g. \cite{BL, GK, Korevaar}. 
\end{itemize}

Similar results, apart from being interesting on their own, have also important applications, as for instance in estimating the spectral gap for the involved operator (see for instance \cite{AC}), 
and for this reason they have been investigated also in Riemannian manifolds (see for instance \cite{SWW} and \cite{Yau}).
On the other hand, geometric structures of Riemannian manifolds are expected to impose strong restrictions on properties of the solutions, and vice versa.
Indeed, Shih~\cite{Shih} has found a bounded convex domain in the hyperbolic plane such that 
the first Dirichlet eigenfunction $\phi$  of $-\Delta$ has a non-convex level set. 
This clearly means that $\phi$ is not  log-concave, and eventually that it is not even  quasi-concave.  
Since quasi-concavity is the weakest among conceivable concavity properties, 
it seems impossible to obtain similar results (as in Euclidean space) about power concavities of solutions 
to elliptic and parabolic boundary value problems on convex domains
in general Riemannian manifolds, especially negatively curved manifolds. 

However, in spite of the Shih result, there are still special situations where some concavity can be expected 
even when the curvature is negative, especially when the domain has a strong symmetry. 
In this paper 
we focus on elliptic and parabolic boundary value problems 
on rotationally symmetric, strongly convex balls in Riemannian manifolds
and obtain strict power concavities of solutions  
even if a Riemannian manifold is negatively curved. 
This may build a foundation for further studying concavity properties of solutions to 
elliptic and parabolic boundary value problems on convex domains in Riemannian manifolds.%, despite the result of Shih.
\medskip

We clarify our setting and notation.
Throughout this paper 
$(M, g)$ is an $N$-dimensional connected, 
complete smooth Riemannian manifold  without boundary, where $N\geq 2$. 
We denote by $d$ the Riemannian distance function on $M$.
Fix $o \in M$ and $R>0$. Set 
\begin{align*}
B(R):=\{x\in M \ |\ d(o,x)<R\},  \quad
\overline{B}(R):=\{x\in M \ |\ d (o,x)\leq R\},
\end{align*}
and 
define the function $\rho:B(R) \to{\bf R}$ by $\rho(x):=d(o,x)$.
Assume that $B(R)$ satisfies the following two conditions: 
\begin{itemize}
  \item[(C1)] 
  $B(R)$ is rotationally symmetric;\\[-19pt]
  \item[(C2)] 
  $B(R)$ is strongly convex.
 \end{itemize}
Under condition~(C1) there exists a unique function $\sigma:(0,R)\to (0,\infty)$ such that 
$B(R)\setminus\{o\}$ is isometric to the warped product $(0,R)\times_\sigma {\bf S}^{N-1}(1)$. 
We call $\sigma$ the \emph{conformal polar factor} of $B(R)$.
A function $u$ in $\overline{B}(R)$ is said to be \emph{rotationally symmetric} 
if there exists a function~$v$ on $[0,R]$ such that $u=v\circ \rho $ on $\overline{B}(R)$. 
\vspace{3pt}

In this paper, under assumptions~(C1) and (C2) 
we obtain sufficient conditions 
for rotationally symmetric solutions to elliptic and parabolic boundary value problems on $B(R)$ 
to be strictly power concave. 
We state the main results of this paper. 
\begin{theorem}
\label{Theorem:1.1}
Under conditions~{\rm (C1)} and {\rm (C2)} 
let $u\in C^2(B(R))\cap C(\overline{B}(R))$ be a rotationally symmetric solution to problem
\begin{equation}
\tag{{\bf E}}
\left\{\begin{array}{ll}
-\Delta  u=F(u)\quad&\mbox{in}\quad B(R),\vspace{3pt}\\
u>0\quad&\mbox{in}\quad  B(R),\vspace{3pt}\\
u=0\quad&\mbox{on}\quad\partial B(R).
\end{array}\right.
\end{equation}
Here $F \in C([0,\infty))\cap C^1((0,\infty))$ is nonnegative on $[0,\infty)$ and $F>0$ on $(0,\infty)$.
Let $\alpha\in[0,1]$ and assume the following conditions.
\begin{itemize}
\item[{\rm (1)}]  
  The function $(0,\infty)\ni s\mapsto s^{\alpha-1}F(s)$ is nonincreasing.
  \item[{\rm (2)}] 
  The conformal polar factor $\sigma$ of $B(R)$ satisfies 
  \begin{equation}
  \label{eq:1.1}
  \frac{d^2}{dr^2}\log \sigma(r)
  \left\{
  \begin{array}{ll}
  \le 0 & \quad\mbox{in}\quad (0,R)\quad\mbox{if}\quad 0 < \alpha<1,\vspace{3pt}\\
  <0 & \quad\mbox{in}\quad (0,R)\quad\mbox{if}\quad \alpha=0,1.
  \end{array}
  \right.
  \end{equation}
  \end{itemize}
Then $u$ is strictly $\alpha$-concave in $B(R)$.
\end{theorem}
\newpage
As a corollary of Theorem~\ref{Theorem:1.1}, we have:  
\begin{corollary}
\label{Corollary:1.1}
Under conditions~{\rm (C1)} and {\rm (C2)} 
let $u\in C^2(B(R))\cap C(\overline{B}(R))$ be a solution to problem
\begin{equation}
\tag{{\bf E}'}
\left\{\begin{array}{ll}
-\Delta  u=\lambda u^\gamma\quad&\mbox{in}\quad B(R),\vspace{5pt}\\
u>0\quad&\mbox{in}\quad  B(R),\vspace{5pt}\\
u=0\quad&\mbox{on}\quad\partial B(R),
\end{array}\right.
\end{equation}
where $\lambda>0$ and $0\le\gamma\le 1$. 
Assume that the conformal polar factor $\sigma$ of $B(R)$ satisfies \eqref{eq:1.1} with $\alpha=1-\gamma$. 
Then $u$ is strictly $(1-\gamma)$-concave in $B(R)$.
\end{corollary}
Problem~({\bf E}') is a generalization of both the Dirichlet eigenvalue problem for $-\Delta$ $(\gamma=1)$ 
and the torsion problem $(\gamma=0)$. 
In the case of $M={\bf R}^N$, we see that $\sigma(r)=r$ and 
$$
\frac{d^2}{dr^2}\log \sigma(r)=-r^{-2}<0\quad\mbox{in}\quad(0,\infty),
$$
so that of course the above results hold. 
%Although we can not say we are giving any new information in this case (since radial solutions in the Euclidean case are well understood and  an explicit representation formula is usually available), we notice that Corollary \ref{Corollary:1.1} for $\gamma<1$ yields a stricter property than the one holding for a generic convex domain (see assertion (a)) and it at least gives an information which does not seem to be always apparent from the explicit representation of a solution, and in some case it may be even sharp in the framework of power concavities (see also Remark~\ref{Remark:1.1} and Appendix).
 
Next we obtain positive power concavity of the first Dirichlet eigenfunction 
under a stronger assumption than \eqref{eq:1.1}. 
\begin{theorem}
\label{Theorem:1.2}
Assume conditions~{\rm (C1)} and {\rm (C2)}. 
If the conformal polar factor $\sigma$  of $B(R)$  satisfies 
\begin{equation}
\label{eq:1.2}
\frac{d^2}{dr^2}\log \sigma(r)\le -\frac{\alpha\lambda_1(B(R))}{N-1}\quad\mbox{in}\quad (0,R)
\end{equation}
for some $\alpha\in(0,1)$, 
then the first Dirichlet eigenfunction  for $-\Delta$ on $B(R)$  is strictly $\alpha$-concave in $B(R)$. 
\end{theorem}
%
%Similarly to the considerations made above for Corollary \ref{Corollary:1.1}, 
We notice that Theorem \ref{Theorem:1.2} may give in Euclidean space an information not immediately retrievable from the explicit representation formula of the first eigenfunction of a ball. See in the Appendix, and in particular Corollary \ref{Corollary:A.1}.

We show that 
a solution to a parabolic boundary value problem with log-concave initial data 
becomes strictly log-concave instantly 
and preserves the strict log-concavity forever.
\begin{theorem}
\label{Theorem:1.3} 
Under conditions~{\rm (C1)} and {\rm (C2)} 
let 
$$
u\in{\mathcal X}:=C^{2,1}(B(R)\times(0,\infty))\cap C(\overline{B}(R)\times(0,\infty))
\cap BC(B(R)\times[0,\infty))
$$
be a nonnegative solution to problem
\begin{equation}
\tag{{\bf P}}
\left\{
\begin{array}{ll}
\partial_t u=\Delta u-G(u) & \mbox{in}\quad B(R)\times(0,\infty),\vspace{3pt}\\
u=0 & \mbox{on}\quad\partial B(R)\times(0,\infty),\vspace{3pt}\\ 
u(x,0)=\varphi(x) & \mbox{in}\quad B(R),
\end{array}
\right.
\end{equation}
where $\varphi$ is a bounded, continuous, nonnegative and rotationally symmetric function in $B(R)$.  
Here $G\in C([0,\infty)]\cap C^2((0,\infty))$. 
Assume the following conditions.
\begin{itemize}
  \item[{\rm (1)}] 
  The conformal polar factor $\sigma$ of $B(R)$ satisfies 
  \begin{equation}
  \label{eq:1.3}
  \frac{d^2}{dr^2}\log \sigma(r)<0\quad\mbox{and}\quad\frac{d^3}{dr^3}\log \sigma(r)\ge 0
  \quad\mbox{in}\quad (0,R).
  \end{equation}
  \item[{\rm (2)}] 
  The function ${\bf R}\ni s\mapsto e^{-s}G(e^s)$ is nonnegative, nondecreasing and convex.
\end{itemize}
Then $u(\cdot,t)$ is strictly log-concave in $B(R)$ for $t>0$ 
provided that $\varphi$ is log-concave in $B(R)$.
\end{theorem}
As a direct consequence of Theorem~\ref{Theorem:1.3}, we have the following result. 
Compare with assertion~(b).  
\begin{corollary}
\label{Corollary:1.2}
Under conditions~{\rm (C1)} and {\rm (C2)} 
let $u\in {\mathcal X}$ be a nonnegative solution to problem
\begin{equation}
\tag{{\bf P}'}
\left\{
\begin{array}{ll}
\partial_t u=\Delta u-\lambda u^\nu & \mbox{in}\quad B(R)\times(0,\infty),\vspace{3pt}\\
u=0 & \mbox{on}\quad\partial B(R)\times(0,\infty),\vspace{3pt}\\ 
u(x,0)=\varphi(x) & \mbox{in}\quad B(R),
\end{array}
\right.
\end{equation}
where $\lambda\ge 0$, $\nu\ge 1$ and $\varphi$ is a bounded, continuous, 
nonnegative and rotationally symmetric function in $B(R)$. 
Assume \eqref{eq:1.3}. 
Then $u(\cdot,t)$ is strictly log-concave in $B(R)$ for $t>0$ 
provided that $\varphi$ is log-concave in $B(R)$.
\end{corollary}
Furthermore, we have the following result for 
an $N$-dimensional simply connected space form  
${\bf M}^N_K$ of nonpositive constant curvature~$K$. 
\begin{corollary}
\label{Corollary:1.3}
Let $\Gamma: {\bf M}_K^N\times{\bf M}_K^N\to {\bf R}$ 
be the heat kernel on ${\bf M}_K^N$, where $K\le 0$.
Then $\Gamma(\cdot,y,t)$ is strictly log-concave on ${\bf M}_K^N$ 
for $y\in{\bf M}_K^N$ and $t>0$. 
\end{corollary}
The heat kernel $\Gamma$ has an explicit representation (see e.g. \cite{GN}). 
Since the representation of $\Gamma$ is complicate for $K<0$, 
it does not seem easy to prove Corollary~\ref{Corollary:1.3} by use of 
direct calculations of second derivatives of $\log\Gamma$. 
\vspace{3pt}

The proofs of our results are based on the standard theory of ODEs and the maximum principle 
for elliptic and parabolic operators. 
We prove Theorems~\ref{Theorem:1.1} and \ref{Theorem:1.2} 
by use of the sign of the first, second and third derivatives of $u^\alpha$ with respect to the radial direction $\rho$. 
For the proof of Theorem~\ref{Theorem:1.3}
we construct approximate solutions $\{u_\epsilon\}$ of a solution~$u$ to problem~({\bf P}) 
and prove that $u_\epsilon(\cdot,t)$ is log-concave in $B(R)$ for $t>0$. 
Then we see that the solution~$u(\cdot,t)$ is log-concave in $B(R)$ for $t>0$. 
Furthermore, we apply the strong maximum principle 
to show that $u(\cdot,t)$ is strictly log-concave in $B(R)$ for $t>0$. 
This completes the proof of Theorem~\ref{Theorem:1.3}. 
\vspace{3pt}

We give some remarks on conditions~(C1), (C2), \eqref{eq:1.1}--\eqref{eq:1.3} 
and concavity properties of the solution to problem~({\bf E}') with $\gamma=0$ on ${\bf R}^N$.
\begin{remark}
\label{Remark:1.1}
\begin{enumerate}
\renewcommand{\labelenumi}{\rm(\roman{enumi})}
\item
Let $(M,g)$ be an $N$-dimensional simply connected space form  
${\bf M}^N_K$ of constant curvature~$K$,  
that is,  
\[
{\bf M}^N_K:=
\begin{cases}
{\bf  S}^N(1/\sqrt{K})  &\text{if $K>0$},\vspace{3pt}\\
{\bf R}^N &\text{if $K=0$},\vspace{3pt}\\
{\bf  H}^N(1/\sqrt{-K}) &\text{if $K<0$}.
\end{cases}
\]
Then conditions~{\rm (C1)}, {\rm (C2)}, \eqref{eq:1.1} and \eqref{eq:1.3} 
are satisfied for balls $B(R)$ if either 
\[
\mbox{{\rm (1)}\quad$K>0$ and $0<R\le\pi/2\sqrt{K}$}
\qquad\mbox{or}\qquad
\mbox{{\rm (2)}\quad$K\le 0$ and $0<R<\infty$}.
\]
For a general Riemannian manifold, 
if $B(R)$ satisfies condition {\rm (C1)},
then there exists $R'\in(0,R]$ such that 
$B(R')$ satisfies {\rm (C1)}, {\rm (C2)}, \eqref{eq:1.1} and \eqref{eq:1.3}. 
On the other hand, condition~\eqref{eq:1.2} is somewhat delicate 
and the constant $\alpha$ in \eqref{eq:1.2} depends on $R$ even for ${\bf M}^N_K$. 
See Appendix for more details.  
\item
Let $M={\bf R}^N$ and  take $o$ as the origin.
Then Corollary~{\rm\ref{Corollary:1.1}} implies that 
the solution to problem~{\rm ({\bf E}')} with $\gamma=0$ is concave. 
This is optimal among power concavities. Indeed, the function $u$ defined by 
$$
u(x)=\frac{\lambda}{2N}(R^2-|x|^2)
$$
is a solution to the  problem and it is $\alpha$-concave in $B(R)$ if and only if $\alpha\le 1$.
\end{enumerate}
\end{remark}

The rest of this paper is organized as follows. 
In Section~2 we recall some properties of rotational symmetry, 
convexities of sets and power concavities for nonnegative functions. 
In Section~3 we prove Theorems~\ref{Theorem:1.1} and \ref{Theorem:1.2}. 
Furthermore, we show that a solution to the parabolic equation corresponding to problem~({\bf E}') 
spontaneously becomes strictly $(1-\gamma)$-concave (see Theorem~\ref{Theorem:3.1}). 
In Section~4 we prove Theorem~\ref{Theorem:1.3}. 
In Appendix we discuss the conditions appearing our results 
in the case of $M={\bf M}_K^N$. 
%
%%%%%%%%%%%%%%%%%%%%%%%%%%%%%%%%%%%%
%%%%%%%%%%%%%%%%%%%%%%%%%%%%%%%%%%%%
\section{Preliminaries}
%%%%%%%%%%%%%%%%%%%%%%%%%%%%%%%%%%%%
%%%%%%%%%%%%%%%%%%%%%%%%%%%%%%%%%%%%
In this section we recall some properties of rotational symmetry
and the notion of some convexities  in Riemannian geometry.

\subsection{Rotational symmetry}
\begin{definition}\label{rotsymdef}
An open metric ball in $M$ centered at $o$ of radius $R$ is said to be \emph{rotationally symmetric}
if there is a function $f:[0,R)\to {\bf R}$ so that $e^{-f(d(o, \cdot ))}g$ is a flat metric on the ball.
\end{definition}
%%%%%
Condition~(C1) guarantees the existence of 
a function $\sigma:(0,R)\to (0,\infty)$ such that $B(R)\setminus\{o\}$ 
is isometric to the warped product $(0,R) \times_\sigma {\bf S}^{N-1}(1)$. 
In this paper we identify  $g$ with
\[
d\rho^2+\sigma(\rho)^2  g_{{\bf S}^{N-1}(1)}
\]
on $B(R)\setminus\{o\}$, 
where $g_{{\bf S}^{N-1}(1)}$ is the standard metric on ${\bf S}^{N-1}(1)$ and $\rho=d(o,\cdot)$. 
Then
  \begin{equation}
  \label{eq:2.1}
  \sigma'(0):=\lim_{r\to +0}\sigma'(r)=1, \quad
  \sigma^{(2n)}(0):=\lim_{r\to +0}\sigma^{(2n)}(r)=0
  \end{equation}
where $n=0,1,2,\dots$ (see \cite[Section 4.3.4]{P}). 
%%%%
%%%
For any smooth rotationally symmetric function $u=v\circ \rho$ 
we have 
\begin{equation}
\label{eq:2.2}
\Delta u=
v''(\rho)+(N-1)(\log\sigma(\rho) )'v'(\rho)
\end{equation}
(see e.g. \cite[p.40, Eq.(29)]{Ch}).
Here and in what follows, $'$ denotes $d/d\rho$. 
%%%%%%%%%%%%%%%%%%%%%%%
\subsection{Convexity}
%%%%%%%%%%%%%%%%%%%%%%%
Let us briefly recall the  notion of convexities in Riemannian geometry.
\begin{definition}\label{Definition:2.2}
\begin{enumerate}
\renewcommand{\labelenumi}{\rm(\roman{enumi})}
\item
We say that a continuous curve $c:[0,1]\to M$ is a \emph{minimal geodesic} if
$$
d(c(s),c(t))=|s-t|\,d(c(0),c(1))\quad\mbox{for}\quad s,t\in[0,1].
$$
%%%%%
\item
We say that a set $\Omega$ of $M$ is \emph{strongly convex}
if for any two points in $\Omega$ there exists a unique minimal geodesic joining them in $M$ 
and the geodesic is contained in $\Omega$.
\item
For any $p\in M$ 
the largest $R$ such that the open metric ball centered at $p$ of radius $R$ 
is strongly convex is called the \emph{convexity radius} at $p$.
\item
Let $\Omega$ be a strongly convex set of $M$.
A function $f:\Omega \to {\bf R}$ is said to be \emph{(strictly) convex} 
if $f \circ c$ is $($strictly$)$ convex for any nonconstant minimal geodesic $c$ in $\Omega$.
\end{enumerate}
\end{definition}
It is worth to mention that the convexity radius is alway positive (see \cite[Theorem IV.5.3]{Sakai}).
\begin{remark} 
\label{Remark:2.1}
Thanks to {\rm\cite[{\it Proposition}~III.~4.8, {\it Lemma}~III.~4.10]{Sakai}}, we see that if
$\Omega$ is a strongly convex set of~$M$ and we fix $p\in \Omega$,
then the following properties hold.
\begin{enumerate}
\renewcommand{\labelenumi}{\rm(\roman{enumi})}
\item
The function $x\mapsto d(p,x)$ is  smooth on $\Omega\setminus\{p\}$ and convex on $\Omega$.
%%%
\item
If a nonconstant minimal geodesic $c:[0,1]\to \Omega$ satisfies  $c(t_0)=p$ for some $t_0\in[0,1]$, 
then $d(p,c(t))=|t-t_0| d(c(0),c(1))$ for $t\in [0,1]$.
\item
If a nonconstant minimal geodesic $c:[0,1]\to \Omega$ satisfies  $p\notin c([0,1])$, 
then 
\[
\frac{d^2}{dt^2}d(p,c(t))>0\quad\mbox{for}\quad t\in[0,1].
\]
\end{enumerate}
\end{remark}
%%%%%%%%%%%%%%%%%%%%%%%
\subsection{Power concavity}
%%%%%%%%%%%%%%%%%%%%%%%
We recall the notion of $\alpha$-concavity for nonnegative functions, 
where $-\infty\le\alpha\le\infty$. 
For $a$, $b>0$, $\mu \in (0,1)$, and $\alpha\in [-\infty, \infty]$,
we define
$$
M_\alpha(a,b;\mu)=\left\{
 \begin{array}{ll}
\left[(1-\mu) a^\alpha+\mu b^\alpha\right]^{1/\alpha}
& \mbox{if}\quad \alpha\not\in\{\pm\infty, 0\},\vspace{5pt}\\
a^{1-\mu}b^{\mu} & \mbox{if}\quad \alpha=0,\vspace{3pt}\\
\max\{a,b\} & \mbox{if}\quad \alpha=\infty,\vspace{3pt}\\
\min\{a,b\}& \mbox{if}\quad \alpha=-\infty,
\end{array}\right.
$$
which is the \emph{$\alpha$\,-\,$($weighted$\,)$ mean }of $a$ and $b$ with ratio $\lambda$. 
For $a$, $b\ge 0$, we define $M_\alpha(a,b;\mu)$ as above if $\alpha\ge 0$ and
$$
M_\alpha(a,b;\mu)=0\quad\mbox{if}\quad \alpha<0\quad\mbox{and}\quad a\cdot b=0.
$$
\begin{definition}
Let $\Omega$ be a strongly convex set of $M$.
Let $u$ be a nonnegative function in $\Omega$ and $-\infty\le\alpha\le\infty$. 
Then we say that $u$ is \emph{$\alpha$\,-concave} 
$($resp.\ \emph{strictly $\alpha$-concave}$)$ in $\Omega$  if 
\[
u(c(t))\ge M_\alpha\big(u(c(0)), u(c(1)); t\big)
\quad
\Big(
\text{resp.\ }
u(c(t)) > M_\alpha\big(u(c(0)), u(c(1)); t\big)
\Big)
\]
for nonconstant minimal geodesics $c:[0,1]\to \Omega$.
In the cases $\alpha=0$ and $\alpha=-\infty$, 
$u$ is also said to be \emph{log-concave} and \emph{quasi-concave} in $\Omega$, respectively. 
\end{definition}
Notice that $\alpha=1$ corresponds to the usual concavity. 
It follows from the Jensen inequality that 
if $v$ is $\alpha$\,-concave in a convex set $\Omega$, 
then $v$ is $\beta$\,-concave in $\Omega$ for any $\beta\le \alpha$.  
This means that quasi-concavity is the weakest among power concavities. 
%%%%%%%%%%%%%%%%%%%%%%%%%%%%%%%%%%%%
%%%%%%%%%%%%%%%%%%%%%%%%%%%%%%%%%%%%
\section{Proofs of Theorems~\ref{Theorem:1.1} and \ref{Theorem:1.2}}
%%%%%%%%%%%%%%%%%%%%%%%%%%%%%%%%%%%%
We prepare key lemmas for the proofs of Theorems~\ref{Theorem:1.1} and \ref{Theorem:1.2}.
%%%%%%%%%%%%%%%%%%%%%%%%%%%%%%%%%%%%
Let $u$ be a rotationally symmetric solution to problem~({\bf E}).
Then there exists $v\in C^2([0,R))\cap C([0,R])$ such that  $u=v\circ \rho$ on $\overline{B}(R)$. 
Furthermore, by \eqref{eq:2.2} we see that $v$ satisfies 
\[
\label{eq:3}
\left\{
\begin{array}{l}
v''+(N-1)(\log\sigma)'v'+F(v)=0\quad\mbox{in}\quad (0,R),\vspace{3pt}\\
v(r)>0\quad\mbox{in}\quad[0,R),\quad v'(0)=v(R)=0.
\end{array}
\right.
\]

Let us introduce the $q$-logarithmic function and the $q$-exponential functions.
We refer to~\cite{Na} for details.
\begin{definition}
For $q\in {\bf R}$, define the \emph{$q$-logarithmic function} on $(0,\infty)$ by
\[ 
L_q(\xi):=\int_1^\xi \frac1{\zeta^q} \,d\zeta. 
\]
The inverse function of $L_q$ is called the \emph{$q$-exponential function}
and it is denoted by $E_q$.
\end{definition}
Note that $L_1$ and $E_1$ correspond to the usual logarithmic function 
and the usual exponential function, respectively.
In the case of $0\leq q<1$ we have
\[ 
L_q(\xi) =\frac{\xi^{1-q}-1}{1-q},\qquad
E_q(\xi) =[1+(1-q)\xi]^{1/(1-q)}: \left(-\frac1{1-q},+\infty\right)\to (0,\infty)\,.
\]
%where $[\xi]_+:=\max\{\xi,0\}$.
For $q\in [0, 1]$, it follows that 
\[
\frac{d}{d\xi} L_q(\xi)=\frac{1}{\xi^q}\ \ \text{for\ }\xi>0,\qquad
\frac{d}{d\xi} E_q(\xi)=E_q(\xi)^q \  \ \text{for\ }\xi >\ell_q:=\begin{cases} -\infty & \text{if}\ q=1, \\-\dfrac{1}{1-q} &\text{if \ $0\leq q<1$}.\end{cases}
\]
Fix $q\in[0,1]$ and set $w:=L_q (v)$.
Then $w$ satisfies
\begin{equation}
\label{eq:3.1}
\left\{
\begin{array}{l}
w''+qE_q(w)^{q-1}( w')^2+(N-1)(\log\sigma)' w'+\dfrac{1}{E_q(w)^q }F(E_q(w))=0 \ \mbox{in}\   (0,R),\vspace{3pt}\\
w(r)>\ell_q,   \quad\mbox{in}\quad[0,R), \quad w'(0)=0, \quad
\lim_{r\to R}w(r)=\ell_q.\vspace{3pt}
\end{array}
\right.
\end{equation}
\begin{lemma}\label{Lemma:3.1}
Assume conditions~{\rm (C1)} and {\rm (C2)}. 
Let $u$ be a rotationally symmetric solution to problem~{\rm({\bf E})} 
and $w$ as above. 
Then $w''(0)<0$ and  $w'<0$ in $(0,R)$. 
\end{lemma}
{\bf Proof.} 
Since $\sigma(0)=0$ and $\sigma'(0)=1$ (see \eqref{eq:2.1}), it follows that 
\begin{equation}
\label{eq:3.2}
(\log \sigma(r))'w'(r)=\sigma'(r)\frac{r}{\sigma(r)-\sigma(0)}\frac{w'(r)-w'(0)}{r}
\to w''(0)\quad\mbox{as}\quad r\to +0.
\end{equation}
By \eqref{eq:3.1} and \eqref{eq:3.2} we have
\[
0=Nw''(0)+\frac{1}{E_q(w(0))^q }F(E_q(w(0)).
\]
Since $F(s)>0$ for $s>0$, 
it follows that $w''(0)<0$. Then we find $\delta\in(0,R)$ 
such that $w'(r)<0$ for $r\in(0,\delta)$. 

By contradiction we prove that $w'<0$ in $(0,R)$.   
Assume that  there exists $r_0\in[\delta,R)$ such that $w'<0$ in $(0,r_0)$ and $w'(r_0)=0$.
Then we have $w''(r_0)\ge 0$. 
On the other hand, we observe from \eqref{eq:3.1} that
\[
w''(r_0)=-\frac{1}{E_q(w(r_0))^q }F(E_q(w(r_0)))<0,
\]
which is a contradiction. 
This means that $w'<0$ in $(0,R)$. 
Thus Lemma~\ref{Lemma:3.1} follows. 
\quad$\Box$
\begin{lemma}\label{Lemma:3.2}
Assume the same conditions as Theorem~{\rm\ref{Theorem:1.1}}. 
Let $u$ be a rotationally symmetric solution to problem~{\rm({\bf E})} 
and $w$ as above with $q=1-\alpha$. 
Then $w''<0$ in $[0,R)$.  
\end{lemma}
{\bf Proof.} 
Assume that the claim  does not hold. 
By Lemma \ref{Lemma:3.1} 
we find $r_1\in(0,R)$ such that 
\begin{equation}
\label{eq:3.3}
w''<0\quad\mbox{in}\quad [0,r_1),\qquad w''(r_1)=0,\qquad w'''(r_1)\ge 0. 
\end{equation}  
On the other hand, 
it follows from \eqref{eq:3.1} that 
\[
\begin{split}
-w'''=&q(q-1)E_q(w)^{2q-1}( w')^3+ 2qE_q(w)^{q-1}w'w''
\\
&+(N-1)(\log\sigma)'' w'
+(N-1)(\log\sigma)' w'' 
-
\frac{qw'}{E_q(w)} F(E_q(w))
+F'(E_q(w))w'
\end{split}
\]
in $(0,R)$. Then, by Lemma \ref{Lemma:3.1} and  \eqref{eq:3.3} 
we have 
\begin{equation}\label{eq:3.4}
0\leq -\frac{w'''}{w'}
=q(q-1)E_q(w)^{2q-1}( w')^2
+(N-1)(\log\sigma)'' 
-\frac{q}{E_q(w)} F(E_q(w))
+F'(E_q(w))
\end{equation}
at $r=r_1$. 
It follows form condition~(2) of Theorem~\ref{Theorem:1.1} that 
\[
-\frac{q}{s} F(s)+F'(s)
=
\frac{1}{s}\left(
(\alpha-1)F(s)+s F'(s)\right)\le 0\quad\mbox{for}\quad s>0. 
\]
Then we deduce from  \eqref{eq:1.1} and \eqref{eq:3.4} that 
\[
0\leq -\frac{w'''}{w'}
\leq q(q-1)E_q(w)^{2q-1}( w')^2
+(N-1)(\log\sigma)'' 
<0\quad\mbox{at}\quad r=r_1,
\]
which is a contradiction. 
Thus Lemma~\ref{Lemma:3.2} follows. 
\quad$\Box$
\vspace{5pt}
\newline
Similarly, we have: 
\begin{lemma}\label{Lemma:3.3}
Assume the same conditions as Theorem~{\rm\ref{Theorem:1.2}}. 
Let $u$ be a rotationally symmetric solution to problem~{\rm({\bf E})} 
with $F(s)=\lambda_1(B(R))s$. 
Let $w$ be as above with $q=1-\alpha$.
Then $w''<0$ in $[0,R)$.  
\end{lemma}
{\bf Proof.}
Similarly to the proof of Lemma~\ref{Lemma:3.2}, 
we assume that the claim does not hold. 
Then we find $r_1\in(0,R)$ such that \eqref{eq:3.3} holds.  
Since $F(s)=\lambda_1(B(R))s$, 
applying a similar argument in the proof of Lemma~\ref{Lemma:3.2}, 
by \eqref{eq:1.2} we obtain
\begin{align*}
0\leq -\frac{w'''}{w'}
&=q(q-1)E_q(w)^{2q-1}( w')^2
+(N-1)(\log\sigma)'' +\alpha  \lambda_1(B(R)) \\
&\leq q(q-1)E_q(w)^{2q-1}( w')^2<0
\end{align*}
at $r=r_1$, instead of \eqref{eq:3.4}. 
This is a contradiction and the proof is complete.
\quad$\Box$
\vspace{4pt}

We complete the proofs of Theorem~\ref{Theorem:1.1}, 
Corollary~\ref{Corollary:1.1} and Theorem~\ref{Theorem:1.2}. 
\vspace{5pt}
\newline
{\bf Proof of Theorem ~\ref{Theorem:1.1}}.  
Let $u$ be a rotationally symmetric solution to problem~({\bf E}). 
Set  $w=L_{1-\alpha} (v)$, where $v$ satisfies $u=v\circ \rho$. 
It follows from Lemmas~\ref{Lemma:3.1} and \ref{Lemma:3.2} that 
$w'<0$ in $(0,R)$ and $w''<0$ in $[0,R)$. 
Let $c:[0,1]\to B(R)$ is a nonconstant minimal geodesic. 
Then, by Remark~\ref{Remark:2.1} 
we observe that 
$$
\frac{d^2}{dt^2}w(\rho(c)))  =
w''(\rho(c))\left(\frac{d}{dt}\rho(c)\right)^2+w'(\rho(c))\frac{d^2}{dt^2}\rho(c)< 0
$$
for $t\in[0,1]$ if $o\notin c([0,1])$ and 
$$
\frac{d^2}{dt^2}w(\rho(c)) =w''(\rho(c))d(c(0),c(1))^2< 0
$$
for $t\in[0,1]\setminus\{t_0\}$ if $o\in c([0,1])$ with $c(t_0)=o$.
These imply that $w(\rho(\cdot))$ is strictly concave in $B(R)$, 
that is, $u$ is strictly $\alpha$-concave in $B(R)$. 
Thus Theorem~\ref{Theorem:1.1} follows.
\quad$\Box$
\vspace{5pt}
\newline
{\bf Proof of Corollary~\ref{Corollary:1.1}}.  
Let $u$ be a solution of problem~({\bf E}) 
with $F(s)=\lambda s^\gamma$, $\lambda>0$ and  $0\le\gamma\leq 1$.
Similarly to \cite[Section~2]{AP} and \cite[Theorem II.5.2]{Ch}, 
we see that $u$ is a unique solution to problem~({\bf E}).
This implies that $u$ is rotationally symmetric in $B(R)$. 
Then we apply Theorem~\ref{Theorem:1.1} with $\alpha=1-\gamma$ 
to see that $u$ is strictly $(1-\gamma)$-concave in $B(R)$. 
Thus Corollary~\ref{Corollary:1.1} follows.
\quad$\Box$
\vspace{5pt}
\newline
{\bf Proof of Theorem~\ref{Theorem:1.2}}. 
By a similar argument to \cite[Theorem II.5.2]{Ch},
we see that the first Dirichlet eigenfunction $u$ for $-\Delta$ on $B(R)$ is rotationally symmetric in $B(R)$.
Let $v$ satisfy $u=v\circ \rho$ and set $w=L_{1-\alpha} (v)$. 
Then we deduce the desired result from Remark~\ref{Remark:2.1} 
together with Lemmas~\ref{Lemma:3.1} and \ref{Lemma:3.2},
similarly to the proof of Theorem~\ref{Theorem:1.1}.\quad$\Box$
\vspace{5pt}

At the end of this section 
we prove the following result, 
in which we combine Corollary~\ref{Corollary:1.1} and the large time behavior 
of a solution to a parabolic boundary value problem.  
\begin{Theorem}
\label{Theorem:3.1}
Under conditions~{\rm (C1)} and {\rm (C2)} 
let $u\in {\mathcal X}$ be a positive solution to 
\begin{equation}
\tag{{\bf P}''}
\left\{
\begin{array}{ll}
\partial_t u=\Delta u+\lambda u^\gamma & \mbox{in}\quad B(R)\times(0,\infty),\vspace{3pt}\\
u=0 & \mbox{on}\quad\partial B(R)\times(0,\infty),\vspace{3pt}\\ 
u(x,0)=\varphi(x) & \mbox{in}\quad B(R),
\end{array}
\right.
\end{equation}
where $\lambda\ge 0$, $0\le\gamma\le 1$ and $\varphi$ 
is a bounded, continuous, nonnegative and rotationally symmetric function in $B(R)$. 
Assume that 
the conformal polar factor $\sigma$ of $B(R)$ satisfies~\eqref{eq:1.1} with $\alpha=1-\gamma$. 
Then there exists $T>0$ such that 
$u(\cdot,t)$ is strictly $(1-\gamma)$-concave in $B(R)$ for $t\ge T$.
\end{Theorem}
{\bf Proof.}
Similarly to \cite[Section~2]{AP}, 
we see that there exists a unique solution $u_\infty\in C^2(\overline{B}(R))$ to problem~({\bf E}'). 
Since $u_\infty$ is rotationally symmetric in $B(R)$, 
there exists $v_\infty\in C^2([0,R])$ such that $u_\infty=v_\infty\circ \rho$ on $\overline{B}(R)$. 
Set $w_\infty:=L_\gamma(v_\infty)$. 
Then, by Lemma~\ref{Lemma:3.2} we see that 
$$
w_\infty''=\frac{v_\infty''v_\infty-\gamma(v'_\infty)^2}{v_\infty^{\gamma+1}}
<0\quad\mbox{on}\quad[0,R).
$$
On the other hand, 
it follows from $u_\infty=0$ on $\partial B(R)$ and the Hopf lemma (see e.g.~\cite[Section~3.2]{GT}) 
that $v_\infty(R)=0$ and $v_\infty'(R)<0$. 
Then we have 
\begin{equation}
\label{eq:3.5}
v_\infty''v_\infty-\gamma(v_\infty')^2<0\quad\mbox{on}\quad[0,R].
\end{equation}

Let $u\in {\mathcal X}$ be a solution to problem~({\bf P}''). 
Let $t_*>0$. Since $v_\infty'(R)<0$, 
we find $m>1$ such that 
$$
0\le u(x,t_*)\le mu_\infty(x)\quad\mbox{on}\quad\overline{B}(R).
$$
Since $u_\infty$ is a solution to problem~({\bf E}'), $\lambda\ge 0$, $0\le\gamma\le 1$ and $m>1$, 
we have
$$
\Delta (mu_\infty)+\lambda(mu_\infty)^\gamma
=\lambda m(-1+m^{\gamma-1})u_\infty^\gamma\le 0\quad\mbox{in}\quad B(R). 
$$
This means that $mu_\infty$ is a supersolution to problem~({\bf P}'') for $t\geq t_*$. 
Let $\underline{u}$ is a positive solution to problem~({\bf P}'') with the zero initial data. 
Then the comparison principle implies that 
\begin{equation}
\label{eq:3.6}
\underline{u}(x,t-t_*)\le u(x,t)\le mu_\infty(x)\quad\mbox{on}\quad\overline{B}(R)\times[t_*,\infty). 
\end{equation}
In particular, $u$ is bounded in $\overline{B}(R)\times[t_*,\infty)$. 
Then, applying the regularity theorems for parabolic equations (see e.g. \cite[Chapter~IV]{LSU}), we see that 
$u$ is smooth in $B(R)\times[t_*,\infty)$ and
\begin{equation}
\label{eq:3.7}
\|u\|_{C^{2,\theta;1,\theta/2}(\overline{B}(R)\times[t_*,\infty))}<\infty
\quad\mbox{for some $\theta\in(0,1)$}.
\end{equation}

Let $\{s_j\}_{j=1}^\infty$ be a sequence on $[t_*,\infty)$ such that $\lim_{j\to\infty}s_j=\infty$ 
and $s_{j+1}>s_j+1$ for $j=1,2,\dots$. 
Set $u_j:=u(x,t+s_j)$. Then $u_j$ satisfies 
\begin{equation}
\label{eq:3.8}
\partial_t u_j=\Delta u_j+\lambda u_j^\gamma\quad\mbox{in}\quad B(R)\times[0,1],
\qquad
u_j=0\quad\mbox{on}\quad \partial B(R)\times[0,1].
\end{equation}
By \eqref{eq:3.7} we apply the Ascoli--Arzel\`a theorem 
to find a subsequence $\{s_{j'}\}$ of $\{s_j\}$ 
and a function $u_*\in C^{2,\theta;1,\theta/2}(\overline{B}(R)\times[0,1])$ such that 
\begin{equation}
\label{eq:3.9}
\lim_{j'\to\infty}\|u_{j'}-u_*\|_{C^{2,\theta;1,\theta/2}(\overline{B}(R)\times[0,1])}=0.
\end{equation}
This together with \eqref{eq:3.8} implies that $u_*$ satisfies 
$$
\partial_t u_*=\Delta u_*+\lambda u_*^\gamma\quad\mbox{in}\quad B(R)\times[0,1],
\qquad
u_*=0\quad\mbox{on}\quad \partial B(R)\times[0,1].
$$
On the other hand, applying the same argument as in the proof of \cite[Theorem~1.1]{INP}, 
we see that $\underline{u}$ converges to $u_\infty$ as $t\to\infty$ uniformly on $\overline{B}(R)$. 
This together with \eqref{eq:3.6} implies that 
\begin{equation}
\label{eq:3.10}
u_\infty(x)\le u_*(x,t)\le mu_\infty(x)\quad\mbox{in}\quad B(R)\times[0,1]. 
\end{equation}
Furthermore, 
since $u$ is smooth in $B(R)\times[t_*,\infty)$, 
we have 
\begin{equation*}
\begin{split}
 & \int_{t_1}^{t_2}\int_{B(R')}(\partial_t u)^2\,dV\,dt
 =\int_{t_1}^{t_2}\int_{B(R')}\partial_t u[\Delta u+\lambda u^\gamma]\,dV\,dt\\
 & =\int_{t_1}^{t_2}\int_{\partial B(R')}\langle\nabla\rho,\nabla u\rangle_g\,\partial_t u\,dS\,dt
 -\int_{t_1}^{t_2}\int_{B(R')}\langle\nabla u,\nabla\partial_t u\rangle_g\,dV\,dt
 +\lambda\int_{t_1}^{t_2}\int_{B(R')}\partial_t u  u^\gamma\,dV\,dt\\
 & =\int_{t_1}^{t_2}\int_{\partial B(R')}\langle\nabla\rho,\nabla u\rangle_g\,\partial_t u\,dS\,dt\\
 & \qquad
 -\frac{1}{2}\int_{t_1}^{t_2}\frac{\partial}{\partial t}\int_{B(R')}|\nabla u|_g^2\,dV\,dt
 +\frac{\lambda}{\gamma+1}\int_{t_1}^{t_2}\frac{\partial}{\partial t}\int_{B(R')}u^{\gamma+1}\,dV\,dt\\
 & \le\int_{t_1}^{t_2}\int_{\partial B(R')}\langle\nabla\rho,\nabla u\rangle_g\,\partial_t u\,dS\,dt
 +\frac{1}{2}\int_{B(R')}|\nabla u(t_1)|_g^2\,dV
 +\frac{\lambda}{1+\gamma}\int_{B(R')}u(t_2)^{1+\gamma}\,dV
\end{split}
\end{equation*}
for $0<R'<R$ and $t_2>t_1\ge t_*$, 
where $dV$ and $dS$ are the volume elements of $(M,g)$ and $\partial B(R')$, respectively. 
Here $\langle\cdot,\cdot\rangle_g$ stands for the Riemannian metric~$g$
and $|\cdot|_g$ is the associated norm. 
Since $\partial_t u=0$ on $\partial B(R)\times(0,\infty)$, 
combining with \eqref{eq:3.7}, we find $C>0$ such that  
\begin{equation*}
\begin{split}
\int_{t_1}^{t_2}\int_{B(R)}(\partial_t u)^2\,dV\,dt
 & \le\frac{1}{2}\int_{B(R)}|\nabla u(t_1)|_g^2\,dV
 +\frac{\lambda}{1+\gamma}\int_{B(R)}u(t_2)^{1+\gamma}\,dV\le C
\end{split}
\end{equation*}
for $t_2\ge t_1\ge t_*$. 
Since $s_{j+1}\ge s_j+1$, we have
$$
\sum_{j'}\int_0^1\int_{B(R)}(\partial_t u_{j'})^2\,dV\,dt
=\sum_{j'}\int_{s_{j'}}^{s_{j'}+1}\int_{B(R)}(\partial_t u)^2\,dV\,dt
\le\int_{t_*}^\infty\int_{B(R)}(\partial_t u)^2\,dV\,dt\le C,
$$
which implies that 
$$
\lim_{j'\to\infty}\int_0^1\int_{B(R)}(\partial_t u_{j'})^2\,dV\,dt=0. 
$$
Then we deduce from \eqref{eq:3.9} that $\partial_t u_*(x,t)=0$ in $B(R)\times[0,1]$, 
that is, $u_*$ is independent of~$t$. 
Furthermore, we observe from \eqref{eq:3.8} and \eqref{eq:3.10} that $u_*$ satisfies
$$
\Delta u_*+\lambda u_*^\gamma=0\quad\mbox{in}\quad B(R),
\qquad 
u_*>0\quad\mbox{in}\quad B(R),
\qquad
u_*=0\quad\mbox{on}\quad\partial B(R). 
$$
The uniqueness of positive solutions to problem~({\bf E}') implies that $u_*=u_\infty$ in $B(R)$. 
Therefore, due to the arbitrariness of $\{s_j\}$, we obtain 
\begin{equation}
\label{eq:3.11}
\lim_{t\to\infty}\|u-u_\infty\|_{C^{2,\theta;1,\theta/2}(\overline{B}(R)\times[t,t+1])}=0.
\end{equation}
By the rotational symmetry of $u$, we find $v\in C^2([0,R]\times(0,\infty))$ such that 
$u(t)=v(t)\circ \rho$ on $\overline{B}(R)\times(0,\infty)$. 
Combining \eqref{eq:3.5} and \eqref{eq:3.11}, 
we find $T>0$ such that 
$$
v''(t)v(t)-\gamma v'(t)^2<0\quad\mbox{on}\quad[0,R]\times[T,\infty). 
$$
Set $w:=L_\gamma(v)$. Then we have 
$$
w''(t)<0\quad\mbox{on}\quad[0,R]\times[T,\infty). 
$$
Similarly to the proof of Theorem~\ref{Theorem:1.1}, 
we observe that $u(t)$ is $(1-\gamma)$-concave in $B(R)$ for $t\ge T$. 
Thus Theorem~\ref{Theorem:3.1} follows.
\quad$\Box$
%%%%%%%%%%%%%%%%%%%%%%%%%%%%%%%%%%%%%%
%%%%%%%%%%%%%%%%%%%%%%%%%%%%%%%%%%%%%%
\section{Proof of Theorem~\ref{Theorem:1.3}}
%%%%%%%%%%%%%%%%%%%%%%%%%%%%%%%%%%%%%%
%%%%%%%%%%%%%%%%%%%%%%%%%%%%%%%%%%%%%%
Assume conditions~(C1) and (C2). 
It follows from condition~(2) of Theorem~\ref{Theorem:1.3} that $G(0)=0$. 
Indeed, if $G(0)>0$, $e^{-s}G(e^s)\to\infty$ as $s\to-\infty$ and $G$ does not satisfy condition~(2). 

Let $\varphi\in BC(B(R))$ be nonnegative, rotationally symmetric and log-concave in $B(R)$.  
Since $G(0)=0$, we can assume, without loss of generality, that $\varphi\not\equiv 0$ in $B(R)$. 
Let $u\in {\mathcal X}$ be a solution to problem~({\bf P}). 
By the uniqueness of solutions to problem~({\bf P}) in $\mathcal X$, 
we see that $u(\cdot,t)$ is rotationally symmetric in $B(R)$ for $t>0$, 
that is, there exists a function 
$$
v\in C^{2,1}([0,R)\times(0,\infty))\cap C([0,R)\times[0,\infty))
$$
such that $u(p,t)=v(\rho(p),t)$ for $p\in \overline{B}(R)\times[0,\infty)$. 
Then $v$ satisfies 
\[
\left\{
\begin{array}{ll}
\partial_t v=v''+(N-1)(\log\sigma)'v'-G(v)\quad & \mbox{in}\quad (0,R)\times(0,\infty),\vspace{3pt}\\
v'(0,t)=v(R,t)=0\quad &\mbox{for}\quad t>0,\vspace{3pt}\\
v(r,0)\ge 0\quad & \mbox{for}\quad r\in [0,R].
\end{array}
\right.
\]
Applying the maximum principle and the Hopf lemma 
(see e.g. \cite[Chapter 2]{Friedman} and \cite[Chapter~II, Section~2]{L}), 
we see that
\begin{equation}
\label{eq:4.1}
v>0\quad\mbox{in}\quad [0,R)\times(0,\infty),
\qquad
v'(R,t)<0\quad\mbox{for}\quad t>0. 
\end{equation} 
Set $w:=\log v$. Then $w$ satisfies 
\begin{equation}
\label{eq:4.2}
\begin{cases}
 \partial_t w
=w''+(w')^2+(N-1)(\log\sigma)'w'-e^{-w}G(e^w)&\mbox{in}\quad(0,R)\times(0,\infty),\vspace{3pt}\\
 w'(0,t)=0, \quad \lim_{r\to R}\,w(r,t)=-\infty&\mbox{for}\quad t>0.\vspace{3pt}
\end{cases}
\end{equation}
We prepare some lemmas for the proof of Theorem~\ref{Theorem:1.3}.
\begin{lemma}
\label{Lemma:4.1}
Assume the same conditions as in Theorem~{\rm\ref{Theorem:1.3}}.  
Further, assume that 
\begin{equation}
\label{eq:4.3}
\varphi\in C^{2,\theta}(\overline{B}(R)),\quad
%\varphi>0\quad\mbox{in}\quad B(R),\quad 
\varphi=0\quad\mbox{and}\quad
\partial_\rho\varphi<0
\quad\mbox{on}\quad\partial B(R),
\end{equation}
where $0<\theta<1$. 
Let $w$ be as above. Then $w'\le0 $ in $[0,R)\times[0,\infty)$.
\end{lemma}
{\bf Proof.} 
Since $\varphi$ is rotationally symmetric and log-concave in $B(R)$, 
we have 
\begin{equation}
\label{eq:4.4}
w''(r,0)\le 0,\qquad w'(r,0)\leq w'(0,0)=0\quad\mbox{for}\quad r\in[0,R]. 
\end{equation}
By \eqref{eq:4.3} we apply the regularity theorems for parabolic equations 
(see \cite[Chapter~IV]{LSU})
to see that 
$$
w\in C^{2,\theta;1,\theta/2}([0,R)\times[0,\infty))\cap C^\infty([0,R)\times(0,\infty)). 
$$
Furthermore, by \eqref{eq:4.2} we have 
\begin{equation}
\label{eq:4.5}
\begin{split}
\partial_t w'=w''' & +2w'w''+(N-1)(\log\sigma)'w''+(N-1)(\log\sigma)''w'\\
& +[e^{-w}G(e^w)-G'(e^w)]w'\quad\mbox{in}\quad(0,R)\times(0,\infty). 
\end{split}
\end{equation}
Let $0<\epsilon<1$ and set $z_1:=w'-\epsilon(1+t)$. 
Since condition~(2) of Theorem~\ref{Theorem:1.3} implies that 
\[
-e^{-s}G(e^s)+G'(e^s)\ge 0\quad\mbox{for}\quad s\in{\bf R}, 
\]
by \eqref{eq:1.3} and \eqref{eq:4.5} we see that $z_1$ satisfies 
\begin{equation}
\label{eq:4.6}
\begin{split}
\partial_t z_1+\epsilon & =z_1''+2w'z_1'+(N-1)(\log\sigma)'z_1'+(N-1)(\log\sigma)''(z_1+\epsilon (1+t))\\
 & \qquad\qquad\qquad\qquad\qquad\quad -[-e^{-w}G(e^w)+G'(e^w)](z_1+\epsilon(1+t))\\
 & < z_1''+2w'z_1'+(N-1)(\log\sigma)'z_1'+(N-1)(\log\sigma)''z_1\\
  & \qquad\qquad\qquad\qquad\qquad\quad 
  -[-e^{-w}G(e^w)+G'(e^w)]z_1
\qquad
\mbox{in}\quad  (0,R)\times(0,\infty).
\end{split}
\end{equation}
Let $0<L<\infty$. 
By \eqref{eq:4.1} and \eqref{eq:4.4}
we find $\delta_1\in(0,R/2)$ such that 
\begin{equation}
\label{eq:4.7}
z_1(r,t)=w'(r,t)-\epsilon(1+t)\le\frac{v'(r,t)}{v(r,t)}-\epsilon<0
\quad
\mbox{in}\quad  \big([0,\delta_1]\cup[R-\delta_1,R)\big)\times[0,L].
\end{equation}
For $\delta\in(0,\delta_1]$, set 
$D_{\delta,L}:=(\delta,R-\delta)\times(0,L]$. 
Assume that 
\[
\max_{(r,t)\in \overline{D_{\delta,L}}}z_1(r,t)\ge 0. 
\]
Then, by \eqref{eq:4.4} and \eqref{eq:4.7}
we find $(r_0,t_0)\in D_{\delta,L}$ such that $z_1<0$ in $[\delta,  R-\delta]\times [0,t_0)$ and $z_1(r_0,t_0)=0$. 
Then 
$$
z_1=z_1'=0,\qquad z_1''\le 0,\qquad\partial_tz_1\ge 0\quad\mbox{at}\quad(r,t)=(r_0,t_0).
$$
These together with \eqref{eq:4.6} imply that 
$$
0<\epsilon\le\partial_t z_1+\epsilon < z_1''\le 0\quad\mbox{at}\quad(r,t)=(r_0,t_0),
$$
which is a contradiction. This means that 
$$
z_1=w'-\epsilon(1+t)<0\quad\mbox{for}\quad (r,t)\in \overline{D_{\delta,L}}=[\delta,R-\delta]\times[0,L].
$$
Since $\delta\in(0,\delta_1)$ is arbitrary, it follows that
$$
z_1=w'-\epsilon(1+t)\le 0\quad\mbox{for}\quad(r,t)\in[0,R)\times[0,L].
$$
Similarly, since $L>0$ and $\epsilon>0$ are arbitrary, 
we see that $w'\le 0$ in $[0,R)\times[0,\infty)$. 
Thus Lemma~\ref{Lemma:4.1} follows.
\quad$\Box$
\begin{lemma}
\label{Lemma:4.2}
Assume the same conditions as in Lemma~{\rm\ref{Lemma:4.1}}. 
Let $\epsilon>0$ and set 
$$
z_2:=w''-\epsilon(1+t).
$$ 
Then, for any $L>0$, 
there exists $\delta_*\in(0,R/2)$ such that 
\[
z_2<0\quad\mbox{in}\quad \left([0,\delta_*)\cup(R-\delta_*,R)\right) \times [0,L]. 
\]
\end{lemma}
{\bf Proof.}
Since $w'(0,t)=0$ for $t\ge 0$ (see \eqref{eq:4.2}), 
it follows from Lemma~\ref{Lemma:4.1} that $w''(0,t)\le 0$ for $t\ge 0$. 
By the continuity of $w''$ we find $\delta_1\in(0,R/2)$ such that  
\begin{equation}
\label{eq:4.8}
z_2(r,t)=w''(r,t)-\epsilon(1+t)<0\quad\mbox{in}\quad 
[0,\delta_1)\times[0,L]. 
\end{equation}
Due to \eqref{eq:4.3} it turns out that $v(R,0)=0$ and $v'(R,0)<0$, 
which implies that there exist  $\delta_2\in(0,\delta_1)$ and $L'\in(0,L)$ such that 
$$
v(r,t)v''(r,t)-v'(r,t)^2<0\quad\mbox{in}\quad [R-\delta_2,R]\times [0,L').
$$
This leads to  
\begin{equation}
\label{eq:4.9}
z_2(r,t)<w''(r,t)=\frac{v(r,t)v''(r,t)-v'(r,t)^2}{v(r,t)^2}<0
\quad
\text{in}\quad  [R-\delta_2,R] \times [0,L'].
\end{equation}

On the other hand,  the boundary condition together with \eqref{eq:4.1}  provides  $v(R,t)=0$ and $v'(R,t)<0$ for $t>0$.
Then, taking a small enough $\delta_3\in(0,\delta_2)$ if necessary, 
we see that 
\[
v(r,t)v''(r,t)-(v'(r,t))^2<0\quad\mbox{in}\quad [R-\delta_3,R] \times [L',L],
\]
which yields
\begin{equation}
\label{eq:4.10}
z_2(r,t)<w''(r,t)=\frac{v(r,t)v''(r,t)-v'(r,t)^2}{z(r,t)^2}<0\quad\mbox{in}\quad [R-\delta_3,R] \times [L',L].
\end{equation}
Combining \eqref{eq:4.8}, \eqref{eq:4.9} and \eqref{eq:4.10},
we complete the proof of Lemma~\ref{Lemma:4.2}.
\quad$\Box$
\begin{lemma}
\label{Lemma:4.3}
Assume the same conditions as in Lemma~{\rm\ref{Lemma:4.1}}. 
Then 
\[
w''\le 0\quad\mbox{in}\quad[0,R)\times[0,\infty). 
\]
\end{lemma}
{\bf Proof.}
The strategy of the proof is as similar as the proof of Lemma \ref{Lemma:4.1}.
Let $L>0$ and $0<\epsilon<2^{-1}(1+L)^{-2}$. 
Set  $z_2:=w''-\epsilon(1+t)$ and take $\delta_*\in(0,R/2)$ as in  Lemma~\ref{Lemma:4.2}.
Since
condition~(ii) of Theorem~\ref{Theorem:1.1} implies that 
$$
-e^{-s}G(e^s)+G'(e^s)\ge  0,\qquad
e^{-s}G(e^s)-G'(e^s)+G''(e^s)e^s\ge 0
\quad
\mbox{for}\quad s\in{\bf R},
$$
by \eqref{eq:1.3}, \eqref{eq:4.5} and Lemma~\ref{Lemma:4.1} 
we have
\begin{align} \notag
 & \epsilon+\partial_t z_2=\partial_t w''\\ \notag
 & =w''''+2(w'')^2+2w'w'''\\ \notag
 & \qquad\quad+(N-1)(\log\sigma)'w'''+2(N-1)(\log\sigma)''w''(N-1)(\log\sigma)'''w'\\ \label{eq:4.11}
 & \qquad\quad-[e^{-w}G(e^w)-G'(e^w)+G''(e^w)e^w](w')^2+[e^{-w}G(e^w)-G'(e^w)]w''\\  \notag 
 & \le z_2''+2(z_2+\epsilon(1+t))^2+2w'z_2'\\ \notag
 & \qquad\quad+(N-1)(\log\sigma)'z_2'+2(N-1)(\log\sigma)''(z_2+\epsilon(1+t))+(N-1)(\log\sigma)'''w'\\  \notag
  & \qquad\quad +[e^{-w}G(e^w)-G'(e^w)](z_2+\epsilon(1+t))\\ \notag
 & \le z_2''+2(z_2+\epsilon(1+t))^2+2w'z_2'+(N-1)(\log\sigma)'z_2'+2(N-1)(\log\sigma)''z_2\\ \notag
 & \qquad\quad+[e^{-w}G(e^w)-G'(e^w)]z_2
\end{align}
for $(r,t)\in(0,R)\times(0,\infty)$.

For $\delta\in (0,\delta_*)$, set $D_{\delta,L}$ be as in the proof of Lemma~\ref{Lemma:4.1}. 
Assume that 
\begin{equation}
\label{eq:4.12}
\max_{(r,t)\in \overline{D_{\delta,L}}}z_2(r,t)\ge 0. 
\end{equation}
By Lemma \ref{Lemma:4.2}, \eqref{eq:4.4} and \eqref{eq:4.12}
we find $(r_0,t_0)\in D_{\delta,L}$ such that $z_2<0$ in $[\delta, R-\delta]\times [0,t_0)$ 
and $z_2(r_0,t_0)=0$. 
Then 
\[
z_2=z_2'=0,\qquad z_2''\le 0,\qquad\partial_tz_2\ge 0\quad\mbox{at}\quad(r,t)=(r_0,t_0).
\]
These together with \eqref{eq:4.11} imply that 
$$
\epsilon\le\epsilon+\partial_t z_2\le z_2''+2\epsilon^2(1+t)^2
\le 2\epsilon^2(1+L)^2<\epsilon\quad\mbox{at}\quad(r,t)=(r_0,t_0),
$$
which is a contradiction. 
This means that 
$$
z_2(r,t)=w''(r,t)-\epsilon(1+t)<0\quad\mbox{for}\quad (r,t)\in\overline{D_{\delta,L}}=[\delta,R-\delta]\times[0,L].
$$
Since $\delta\in(0,\delta_*)$ is arbitrary, we have
$$
z_2(r,t)=w''(r,t)-\epsilon(1+t)\le 0\quad\mbox{in}\quad  [0,R)\times[0,L].
$$
Letting $\epsilon\to +0$, we see that
$$
w''(r,t)\le 0\quad\mbox{in}\quad  [0,R)\times[0,L].
$$
Since $L$ is arbitrary, 
we obtain the desired conclusion. 
Thus the proof is complete. 
\quad$\Box$
\vspace{5pt}

Now we are ready to complete the proof of Theorem~\ref{Theorem:1.3}.
\vspace{5pt}
\newline
{\bf Proof of Theorem~\ref{Theorem:1.3}}. 
Let $\varphi\in BC(B(R))$ be nonnegative, rotationally symmetric and log-concave in $B(R)$. 
We can assume, without loss of generality, that $\varphi\not\equiv 0$ in $B(R)$. 
There exists $\psi:[0,R) \to {\bf R}$ such that  $\varphi=\psi \circ\rho$ on  $B(R)$. 
Let $u$ be a solution to problem~({\bf P}) and set $w:=\log v$, where $v$ is as above. 

Let  $B_R$ denote the open metric ball of radius $R$, centered  at the origin in ${\bf R}^N$.
Let  $\eta$ be a solution of 
\[
\left\{
\begin{array}{ll}
\partial_t\eta=\Delta\eta & \quad\mbox{in}\quad B_R\times(0,\infty),\vspace{3pt}\\
\eta=0 & \quad\mbox{on}\quad \partial B_R\times(0,\infty),\vspace{3pt}\\
\eta(x,0)=\psi(|x|)& \quad\mbox{in}\quad \partial B_R.
\end{array}
\right.
\]
Since $\eta(\cdot, 0)$ is log-concave in $B_R\subset {\bf R}^N$, 
by assertion~(b) in Section~1 
we see that $\eta(\cdot,t)$ is log-concave in $B_R$ for $t>0$.
Furthermore, 
we deduce from the Hopf lemma and the regularity theorems for parabolic equations that
\begin{equation}
\label{eq:4.13}
\eta(t)\in C^{2,\theta}(\overline{B_R}),\qquad 
\partial_r\eta(t)<0\quad\mbox{on}\quad\partial B_R,
\end{equation}
for $t>0$, where $0<\theta<1$. 

Let $\epsilon>0$. 
Condition~(ii) of Theorem~\ref{Theorem:1.1} implies that $G'\ge 0$ in $(0,\infty)$. 
Then, similarly to \cite{FH}, 
we find a unique solution $u_\epsilon\in{\mathcal X}$ 
to problem~({\bf P}) with $\varphi$ replaced by $\eta(\rho, \epsilon)$. 
The uniqueness of solutions to problem~({\bf P}) yields that 
$u_\epsilon(\cdot,t)$ is rotationally symmetric in $B(R)$ for $t\ge 0$. 
This implies that there exists a function $v_\epsilon$ in $[0,R]\times[0,\infty)$ 
such that $u_\epsilon(p,t)=v_\epsilon(\rho(p),t)$ 
for $p\in \overline{B}(R)$ and $t\ge 0$. 
Set $w_\epsilon:=\log v_\epsilon$. 
By \eqref{eq:4.13} we apply Lemmas~\ref{Lemma:4.1} and \ref{Lemma:4.3} to obtain 
$$
w_\epsilon'\le 0\quad\mbox{and}\quad
w_\epsilon''\le 0\quad\mbox{in}\quad [0,R)\times[0,\infty). 
$$
Since $v_\epsilon(t)\to v(t)$ as $\epsilon\to 0$ in $C^2([0,R])$ 
for $t>0$, by \eqref{eq:4.4} we have 
\begin{equation}
\label{eq:4.14}
w'\le 0\quad\mbox{and}\quad w''\le 0\quad\mbox{in}\quad [0,R)\times[0,\infty). 
\end{equation}

We prove that 
\begin{equation}
\label{eq:4.15}
w''<0\quad\mbox{in}\quad (0,R)\times[0,\infty). 
\end{equation}
Set $z_2=w''$. 
Similarly to \eqref{eq:4.11}, 
by \eqref{eq:1.3} and \eqref{eq:4.14} it turns out that 
$$
\partial_t z_2
\le z_2''+(N-1)(\log\sigma)'z_2'+2(N-1)(\log\sigma)''z_2+2z_2^2+2w'z_2'+[e^{-w}G(e^w)-G'(e^w)]z_2
$$
for $(r,t)\in(0,R)\times(0,\infty)$. 
Furthermore, it follows from \eqref{eq:4.14} that $z_2\le 0$ in $(0,R)\times(0,\infty)$. 
Since $z_2\not\equiv 0$ in $(0,R)\times(0,\infty)$, 
we see that $z_2\not\equiv 0$ in $(\delta,R)\times(0,\infty)$ 
for small enough $\delta\in(0,R)$. 
Then, applying the strong maximum principle (see e.g. \cite[Chapter~2, Section~2, Theorem~3 ]{Friedman}),  
we see that $z_2<0$ in $(\delta,R)\times(0,\infty)$. 
Since $\delta$ is arbitrary, we have 
$w''=z_2<0$ in $(0,R)\times(0,\infty)$, which implies \eqref{eq:4.15}. 
Then, by the same argument of the proof of Theorem~\ref{Theorem:1.1} 
we see that $u(\cdot,t)$ is strictly log-concave in $B(R)$ for $t>0$. 
Thus Theorem~\ref{Theorem:1.3} follows.
\quad$\Box$
\vspace{5pt}%

Combining Theorem~\ref{Theorem:1.2}, Theorem~\ref{Theorem:3.1}
and Corollary~\ref{Corollary:1.2}, 
we have:
\begin{corollary}
\label{Corollary:4.1}
Under conditions~{\rm (C1)} and {\rm (C2)} 
let $u\in {\mathcal X}$ be a nonnegative solution 
to the problem
\begin{equation*}
\left\{
\begin{array}{ll}
\partial_t u=\Delta u & \mbox{in}\quad B(R)\times(0,\infty),\vspace{3pt}\\
u=0 & \mbox{on}\quad\partial B(R)\times(0,\infty),\vspace{3pt}\\ 
u(x,0)=\varphi(x) & \mbox{in}\quad B(R),
\end{array}
\right.
\end{equation*}
where $\varphi$ is a bounded continuous, nonnegative and rotationally symmetric function in $B(R)$. 
Then the following properties hold. 
\begin{itemize}
  \item[{\rm (i)}] 
  Assume \eqref{eq:1.1} with $\alpha=0$. 
  Then there exists $T_1>0$ such that $u(\cdot,t)$ is strictly log-concave in $B(R)$ for $t\ge T_1$.
  \item[{\rm (ii)}] 
  Assume \eqref{eq:1.2} for some $\alpha\in(0,1)$. 
  Then there exists $T_2>0$ such that $u(\cdot,t)$ is strictly $\alpha$-concave in $B(R)$ for $t\ge T_2$. 
  \item[{\rm (iii)}] 
  Assume \eqref{eq:1.3}. 
  Then  $u(\cdot,t)$ is strictly log-concave in $B(R)$ for $t>0$ 
  provided that $\varphi$ is rotationally symmetric and log-concave in $B(R)$. 
\end{itemize}
\end{corollary}
{\bf Proof}. 
Assertion~(iii) immediately follows from Theorem~\ref{Theorem:1.3}. 
Let $\tilde{u}:=e^{\lambda_1(B(R))t}u$. Then $\tilde{u}$ satisfies 
\begin{equation*}
\left\{
\begin{array}{ll}
\partial_t \tilde{u}=\Delta \tilde{u}+\lambda_1(B(R))\tilde{u} & \mbox{in}\quad B(R)\times(0,\infty),\vspace{3pt}\\
\tilde{u}=0 & \mbox{on}\quad\partial B(R)\times(0,\infty),\vspace{3pt}\\ 
\tilde{u}(x,0)=\varphi(x)\ge 0 & \mbox{in}\quad B(R). 
\end{array}
\right.
\end{equation*}
Assume \eqref{eq:1.1} with $\alpha=0$.
By Theorem~\ref{Theorem:3.1} with $\gamma=1$ 
we find $T_1>0$ such that $\tilde{u}(t)$ is log-concave in $B(R)$ for $t\ge T_1$. 
This implies assertion~(i). 

It remains to prove assertion~(ii).
Assume \eqref{eq:1.2} for some $\alpha\in(0,1)$. 
Using Lemma~\ref{Lemma:3.3}, instead of Lemma~\ref{Lemma:3.2}, 
and applying a similar argument to that of the proof of Theorem~\ref{Theorem:3.1}, 
we find $T_2>0$ such that $\tilde{u}(t)$ is $\alpha$-concave in $B(R)$ for $t\ge T_2$. 
This implies assertion~(ii). Thus Corollary~\ref{Corollary:4.1} follows. 
\quad$\Box$\vspace{5pt}

Finally we prove Corollary~\ref{Corollary:1.3}.
\vspace{5pt}\newline
{\bf Proof of Corollary~\ref{Corollary:1.3}.}
For any $\epsilon>0$
let $\varphi_\epsilon$ be a continuous, nonnegative, rotationally symmetric and log-concave function in ${\bf M}_K^N$ 
such that 
\begin{equation}
\label{eq:4.16}
\mbox{supp}\,\varphi_\epsilon\subset B(\epsilon),\qquad \int_{{\bf M}_K^N}\varphi_\epsilon(x)\,dV=1.
\end{equation}
For any $n=1,2,\dots$, let $u_{n,\epsilon}$ be a solution to the problem
\begin{equation*}
\left\{
\begin{array}{ll}
\partial_t u=\Delta u & \quad\mbox{in}\quad B(n)\times(0,\infty),\vspace{3pt}\\
u=0 & \quad\mbox{on}\quad\partial B(n)\times(0,\infty),\vspace{3pt}\\
u(x,0)=\varphi_\epsilon(x) & \quad\mbox{in}\quad B(n).
\end{array}
\right.
\end{equation*}
Applying the standard arguments for parabolic equations, 
we see that 
$$
\lim_{n\to\infty}u_{n,\epsilon}(x,t)
=\int_{{\bf M}_K^N}\Gamma(x,y,t)\varphi_\epsilon(y)\,dV
$$
for $x\in{\bf M}_K^N$ and $t>0$.
Then we deduce from \eqref{eq:4.16} that 
\begin{equation}
\label{eq:4.17}
\lim_{\epsilon\to+0}\lim_{n\to\infty}u_{n,\epsilon}(x,t)=\Gamma(x,o,t)
\end{equation}
for $x\in{\bf M}_K^N$ and $t>0$. 
Since $u_{n,\epsilon}(\cdot,t)$ is log-concave in $B(n)$ for $t>0$, $n=1,2,\dots$ and $\epsilon>0$ 
(see Corollary~\ref{Corollary:4.1}~(iii)), 
we observe from \eqref{eq:4.17} that $\Gamma(\cdot,o,t)$ is log-concave in ${\bf M}_K^N$ for $t>0$. 
Furthermore, by the same argument as in the proof of \eqref{eq:4.15} 
we see that $\Gamma(\cdot,o,t)$ is strictly log-concave in ${\bf M}_K^N$ for $t>0$. 
Then, combining the arbitrariness of $o\in{\bf M}_K^N$, 
we complete the proof of Corollary~\ref{Corollary:1.3}. 
\quad$\Box$
%%%%%%%%%%%%%%%%%%%%%%%%%%%%%%%%%%%%%%%%%%%%%%%%%%%%%%%%%
%%%%%%%%%%%%%%%%%%%%%%%%%%%%%%%%%%%%%%%%%%%%%%%%%%%%%%%%%
\renewcommand{\thesection}{\Alph{section}}
\setcounter{section}{0}
\section{Appendix}
%%%%%%%%%%%%%%%%%%%%%%%%%%%%%%%%%%%%%%%%%%%%%%%%%%%%%%%%%
%%%%%%%%%%%%%%%%%%%%%%%%%%%%%%%%%%%%%%%%%%%%%%%%%%%%%%%%%
We  discuss the conditions (C1), (C2),  \eqref{eq:1.1},  \eqref{eq:1.2} and \eqref{eq:1.3} in our results .
Note that if $B(R)$ satisfies (C1), 
then $B(R')$ satisfies (C1) and (C2) for small enough $R'\in(0,R)$ 
since the convexity radius never vanishes.
Since 
\begin{equation*}
\begin{split}
(\log \sigma)'(r)&=\frac{\sigma'(r)}{\sigma(r)},\qquad
(\log \sigma)''(r)=\frac{\sigma(r)\sigma''(r)-\sigma'(r)^2}{\sigma(r)^2}, \\
(\log \sigma)'''(r)&=\frac{2\sigma'(r)^3  -3\sigma(r)\sigma'(r)\sigma''(r) +\sigma^2(r)\sigma'''(r)  }{\sigma(r)^3},
\end{split}
\end{equation*}
it follows from \eqref{eq:2.1} 
that \eqref{eq:1.1} and \eqref{eq:1.3} hold for small enough $R'\in(0,R]$.

Let us first recall some properties of geodesics and sectional curvatures in rotationally symmetric metrics.
%%%
\begin{proposition}\,{\rm (c.f. \cite[Propositions~7.38, 7.42]{Oni})}\label{Proposition:A.1}\!
Under the identification of $B(R)\setminus\{o\}$ with the warped product  $(0,R) \times_\sigma{\bf S}^{N-1}(1)$,
a curve $(a,b)$ in $(0,R) \times_\sigma{\bf S}^{N-1}(1)$ is a geodesic if and only if 
\begin{equation}\label{eq:A.2}
a''=|b'|^2_{{\bf S}^{N-1}} \sigma(a) \sigma'(a)
\quad\mbox{and}\quad
b''=\frac{-2}{\sigma(a)}\sigma'(a)a'b'.
\end{equation}
%%%%%%%%%%%
Here $|\cdot|_{{\bf S}^{N-1}}$ is the norm associated to $g_{{\bf S}^{N-1}(1)}$. 
The sectional curvature $K_{(r,\xi)}$ of the tangent plane 
containing $\partial_r$ at $(r,\xi) \in B(R)\setminus\{o\}$ is given by
\[
  K_{(r,\xi)}=
  -\frac{\sigma''(r)}{\sigma(r)}.
\]
\end{proposition}
%%%%%%%%%%%%%%%%%%%%%%%%%
We observe from \eqref{eq:2.1} that 
\[
\lim_{r \to +0}  K_{(r,\xi)}=- \lim_{r \to +0} \sigma'''(r).
\]
Regarding $B(R)$ as $[0,R) \times {\bf S}^{N-1}(1)$, 
for any $v\in {\bf S}^{N-1}(1)$ and $0\leq R_\pm <R$ 
we easily check that 
\[
[0,R_+) \times \{ v \} \cup  [0, R_-) \times \{  -v \} 
\]
is the imagine of a geodesic of length $(R_++R_-)$.

Under condition  {\rm (C1)} let us consider the relation 
between the strong convexity of $B(R)$ and  the positivity of $\sigma'$ on $(0,R)$.
%%%%%%%%%%%%%%%%%%%%%%%%%%%%%% 
\begin{lemma}\label{Lemma:A.1}
Assume condition~{\rm (C1)}. 
If $B(R)$ is strongly convex, then $\sigma'>0$ on~$(0,R)$.
Conversely, if $R=\infty$ and  $\sigma'>0$ on~$(0,\infty)$, 
then $B(R')$ is strongly convex for $R'>0$.
\end{lemma}
{\bf Proof.}
Let $B(R)$ be strongly convex. Assume that there exists $R_0\in (0,R)$ such that
\[
\sigma'(R_0)=0, \qquad
\sigma'(r)>0\quad \text{in}\ (0,R_0).
\]
Let us regard $B(R)$ as  $[0,R) \times {\bf S}^{N-1}(1)$ and take $p_{\pm}=(R_0, \pm e_1) \in  B(R)$.
If a curve $(a_1,b_1): I \to B(R)$ joins  $p_+$ and $p_-$, then so does the curve $(a_1,Pb_1): I \to B(R)$ for $P\in O(N)$ such that $Pe_1=e_1$.
Since the curves have the same length, 
the uniqueness of minimal geodesics implies that 
the image of the  minimal geodesic joining $p_+$ and $p_-$ is  
\[
\big([0,R_0) \times \{ e_1 \} \big)\cup \big( [0, R_0) \times \{  -e_1 \}\big) ,
\]
consequently $d(p_+, p_-)=2R_0$.
On the other hand,  a curve $c=(a,b): [0,\pi] \to B(R)\setminus\{o\} $ of the form 
\[
a(t) \equiv R_0, \qquad
b(t)=
e_1 \cos t+e_2 \sin t 
\]
is a  geodesic from $p_+$ to $p_-$ by  Proposition \ref{eq:A.2}.
By the uniqueness of minimal geodesics joining $p_+$ and $p_-$
we have
\[
2R_0=d(p_+, p_-)<\int_0^\pi \sqrt{g(c'(t), c'(t))} dt = \int_0^{\pi} s(R_0) dt=\pi s(R_0).
\]
Let us set
\[
\theta:=\frac{2R_0}{s(R_0)}<\pi, 
\quad p:=c(\theta),
\]
and denote by $(a_2,b_2): I \to B(R)$ a minimal geodesic from  $p_+$ and $p$.
For $P_2\in O(N)$ which fixes $e_1, e_2$, we observe that $(a_2,P_2b_2): I \to B(R)$ is also minimal geodesic from  $p_+$ and $p$.
This implies that the restriction of $c$ to $[0,\theta]$ is a minimal geodesic from  $p_+$ to $p$ and $d(p_+,p)=2R_0$.
On the other hand,  a curve  $c_3=(a_3,b_3): [-R_0, R_0] \to B(R)$ of the form 
\[
a_3(t)=|t|, \quad
b_3(t)=\begin{cases} e_1, & -R_0\leq t<0 \\ b(\theta) ,&0 \leq t\leq R_0\end{cases}
\]
is a broken geodesic from $p$ to $p_+$.
%%%%%
Since a broken geodesic is not minimal, it turns out that 
 \begin{align*}
2R_0=d(p_+, p)
&<
\int_{-R_0}^0 \sqrt{g(c_3'(t),c_3'(t))}   dt
+
\int_{-R_0}^0 \sqrt{g(c_3'(t),c_3'(t))}   dt
=2R_0,
\end{align*}
%%%%
which is a contradiction.
Thus we conclude $\sigma'>0$  on $(0,R)$.

Conversely, assume $R=\infty$  and  the positivity of $\sigma'$ on  $(0,\infty)$.
Let $c:[0,1] \to M$ is a minimal geodesic with $c(0)$, $c(1)\in B(R')$, where $R'>0$. 
Then we derive from the positivity of  $\sigma'$ and \eqref{eq:A.2} 
for the radial part together  with the minimality of $c$  that $c([0,1])\subset B(R')$
and the uniqueness of a unique minimal geodesic joining $c(0)$ and $c(1)$.
Thus $B(R')$ is strong convex. 
\quad$\Box$\vspace{5pt}

We recall the comparison theorem  for the first eigenvalue on $B(R)$. 
Let $K_{\max}(R)$ (resp. $K_{\min}(R)$) be the maximum (resp. minimum) of 
the sectional curvature of the tangent plane containing $\partial_\rho$ on $\overline{B}(R)$. 
%%%%
\begin{proposition} {\rm (\cite[Theorem~3.6]{Cheng1})}\label{Proposition:A.2}
Let $\lambda_1(K,R,N)$ denote the first Dirichlet eigenvalue of the open metric ball of radius $R$ in ${\bf M}^N_K$.
Then 
\[
\lambda_1(B(R)) \geq  \lambda_1(K_{\max}(R), R,N)\,.
\]
\end{proposition}
Assume that $B(R)$ satisfies (C1) and (C2).
Then we observe from (A1) that 
$\alpha$ satisfies \eqref{eq:1.2} if and only if 
\begin{align}
\begin{split}
\label{eq:A.3}
0<\alpha\le  A(\sigma,R,N)
:&= -\frac{N-1}{\lambda_1(B(R))}\sup_{0< r< R}(\log\sigma)'' \\ 
&=\frac{N-1}{\lambda_1(B(R))}\inf_{0<r< R}\left(-\frac{\sigma''(r)}{\sigma(r)}-\frac{1-\sigma'(r)^2}{\sigma(r)^2}+\frac{1}{\sigma(r)^2}\right).
\end{split}
\end{align}
It follows that
\[
K_{\min}(R) \leq -\frac{\sigma''(r)}{\sigma(r)} \leq K_{\max}(R)
\quad
\text{on}\ [0,R].
\]
By condition~(C1) and Lemma \ref{Lemma:A.1} 
we see that $\sigma$ and $\sigma'$ are positive on $(0,R)$. 
This implies that
\[
 K_{\min}(R) \sigma(r)\sigma'(r) \leq -\sigma'(r)\sigma''(r) \leq K_{\max}(R)\sigma(r)\sigma'(r)
\quad
\text{on}\ [0,R].
\]
Integrating the above inequality on $[0,r]$ with condition \eqref{eq:2.1} 
and dividing by $\sigma(r)^2$ provide
\[
 K_{\min}(R)  \leq  \frac{1-\sigma'(r)^2}{\sigma(r)^2}  \leq K_{\max}(R)
\quad
\text{on}\ [0,R].
\]
Therefore, by \eqref{eq:A.3} we obtain
\begin{align*}
A(\sigma,R,N)
&\leq \frac{N-1}{\lambda_1(B(R))}\left(K_{\max}(R)-K_{\min}(R)+\frac{1}{\sigma(R)^2}\right) \\
&\leq  \frac{N-1}{\lambda_1(K_{\max}(R),R,N)} \left(K_{\max}(R)-K_{\min}(R)+\frac{1}{\sigma(R)^2}\right).
\end{align*}
Thus $A(\sigma,R,N)$ can be estimated in terms of the curvature bounds.

Let $M$ be an $N$-dimensional simply connected space form  ${\bf M}^N_K$ of constant curvature $K$,  
that is, ${\bf M}^N_K$ is as in Remark~\ref{Remark:1.1}~(i). 
The convexity radius  of ${\bf M}^N_K$ is independent of the choice of points, which is given by
\[
r_K:=\begin{cases}\dfrac{\pi}{2\sqrt{K}}&\text{if $K>0$},\\ \infty&\text{if $K\leq 0$.}\end{cases}
\]
The conformal polar factor $\sigma$ in ${\bf M}^N_K$ is given by 
\[
\sigma(r)=\sigma_K(r):=\begin{cases}
\dfrac{\sin(\sqrt{K}r)}{\sqrt{K}} & \text{if $K>0$},\vspace{3pt}\\
 r & \text{if $K= 0 $},\vspace{3pt}\\
\dfrac{\sinh(\sqrt{-K}r)}{\sqrt{-K}} & \text{if $K<0$}.
\end{cases}
\]
%%%%%%%%
This implies that 
$$
\sigma_K'(r)=
\left\{
\begin{array}{ll}
\cos \sqrt{K}r & \mbox{if $K>0$},\vspace{3pt}\\
1 & \mbox{if $K=0$},\vspace{3pt}\\
\cosh \sqrt{-K}r & \mbox{if $K<0$},\vspace{3pt}\\
\end{array}
\right.
$$
hence $\sigma_K'$ is positive in $(0,r_K)$.
%%%%%
Then it turns out that 
\[
  (\log\sigma_K)'=\frac{\sigma_K'}{\sigma_K}>0,\quad
 (\log\sigma_K)''
 =-\frac{1}{\sigma_K^2}<0, \quad
  (\log\sigma_K)'''=\frac{2\sigma_K'}{\sigma_K^{3}}>0
\text{\quad in\ }(0,r_K),
\]
since  
$\sigma_K''=-K\sigma_K$ and 
$K\sigma_K^2+(\sigma_K'')^2=1$. 
%%%%%
Therefore
all the conditions (C1), (C2)
\eqref{eq:1.1} and \eqref{eq:1.3} for $B(R)$ in ${\bf M}_K^N$ are satisfied if $R\in{\bf R}$ with $0<R\leq r_K$. 
Furthermore, we observe that 
\[
A(\sigma_K,R,N)
= \frac{N-1}{\lambda_1(K,R,N) \sigma_K(R)^2}
=\left\{
\begin{array}{ll}
\displaystyle{\frac{K(N-1)}{\lambda_1(K,R,N)\sin^2(\sqrt{K}R )}}\quad & \mbox{if}\quad K>0,\vspace{7pt}\\
\displaystyle{\frac{N-1}{\lambda_1(0,R,N)R^2}} & \mbox{if}\quad K=0,\vspace{7pt}\\
\displaystyle{\frac{-K(N-1)}{\lambda_1(K,R,N)\sinh^2(\sqrt{-K}R )}}\quad & \mbox{if}\quad K<0.\vspace{3pt}
\end{array}
\right.
\]
In particular, in the case of $K=0$, it follows from \cite[Theorem II.5.4]{Ch} that 
$$
A(\sigma_0, R, N)=\frac{N-1}{\lambda_1(0,R,N)R^{2}}
=\frac{N-1}{j_{(N-2)/2}^2},
$$
which is independent of $R$. 
Here $j_{(N-2)/2}$ is the first positive zero of the Bessel function $J_a$ of the first kind, where $a=(N-2)/2$. 
Consequently, we have:
\begin{corollary}
\label{Corollary:A.1}
Let $N\ge 2$ and $R$, $K\in{\bf R}$ with $0<R\le r_K$. 
Let $\alpha\in(0,1)$ be such that $\alpha\le A(\sigma_K,R,N)$. 
Then the first Dirichlet eigenfunction for $-\Delta$ on $B(R)$ in ${\bf M}_K^N$ 
is strictly $\alpha$-concave in $B(R)$. 
\end{corollary}
Although the first Dirichlet eigenfunction on a convex domain in the hyperbolic plane 
is not necessarily quasi-concave (see \cite{Shih}), 
Corollary~\ref{Corollary:A.1} says that 
the first Dirichlet eigenfunction for $-\Delta$ on $B(R)$ of ${\bf M}_K^N$, where $0<R<\infty$, 
is positive power concave for $0<R\le r_K$ even if $K<0$. 

On general Riemannian manifolds, 
since 
$$
\lim_{R\to 0+}R^2\cdot \lambda_1(B(R))=\lambda_1(0, 1, N)
$$
(see \cite[Problem for Chapter~VI.9]{Sakai}), 
by \eqref{eq:2.1} we see that 
\begin{align*}
\lim_{R\to 0+}A(\sigma,R,N)
&=\lim_{R\to 0+} \left\{\frac{N-1}{R^2 \lambda_1(B(R))} \inf_{0<r<R} R^2 \cdot \left(-\frac{\sigma''(r)}{\sigma(r)}-\frac{1-\sigma'(r)^2}{\sigma(r)^2}+\frac{1}{\sigma(r)^2}\right)\right\}\\
&=\frac{N-1}{\lambda_1(0,1,N)}>0.
\end{align*}
Consequently we have: 
\begin{corollary}
\label{Corollary:A.2}
Let $N\ge 2$. Assume that $B(R)$ satisfies condition~{\rm (C1)}. 
Let $\alpha\in(0,1)$ be such that 
$$
\alpha<\frac{N-1}{\lambda_1(0,1,N)}. 
$$
Then there exists $R_*>0$ such that 
the first Dirichlet eigenfunction for $-\Delta$ on $B(R')$ 
is strictly $\alpha$-concave in $B(R')$ for $0<R'\le R_*$.
\end{corollary}
%%%%%%%%%%%%%%%%%%%%%%%%%%%%%%%%%%%%%%%%%%%%%%%%%
\vspace{3pt}
%\newline
%%%%%%%%%%%%%%%%%%%%%%%%%%%%%%%%%%%%%%%%%%%%%%%%
{\bf Acknowledgements.} 
The first author was supported in part by the Grant-in-Aid for Scientific Research (S)(No.~19H05599)
from Japan Society for the Promotion of Science.
The second author has been partially supported by INdAM through a GNAMPA Project.
The third author was supported in part 
by the Grant-in-Aid for Scientific Research (C)(No.~19K03494) 
and
the Grant-in-Aid for Scientific Research (C)(No.~19H01800).
%%%%%%%%%%%%%%%%%%%%%%%%%%%%%%%%%%%%%%
%%%%%%%%%%%%    references    %%%%%%%%%%%%%%%%%%
%%%%%%%%%%%%%%%%%%%%%%%%%%%%%%%%%%%%%%

%%%%%%%%%%%%%%%%%%%%%%%%%%%%%%%%%%%%
%%%%%%%%%%%%%%%%%%%%%%%%%%%%%%%%%%%%
\end{document}